\documentclass[11pt]{article}
\usepackage{cite}
\usepackage{amsmath}
\usepackage{amssymb}
\usepackage{mathtools}
\usepackage{graphicx}
\usepackage{color}
\usepackage{enumitem}
\usepackage{tabularx}
\usepackage{listings}

\usepackage{algorithm}
\usepackage{placeins} 
\usepackage{array} 
\usepackage{algpseudocode}
\algnewcommand{\LineComment}[1]{\State \(\triangleright\) {#1}}

\usepackage{float}
\newfloat{algorithm}{t}{lop}

\usepackage[top=1.5in,bottom=1in,right=1in,left=1in]{geometry}
\usepackage[hidelinks]{hyperref}
\usepackage{xspace}
\usepackage{subcaption}

\usepackage[nottoc]{tocbibind}

\usepackage{todonotes}
\usepackage[final]{showlabels}

\newcommand{\R}{\mathbb{R}}

\newcommand{\Si}{\mathbb{S}_{i}^{i+1}}

\newcommand{\hi}[1]{h^{(#1)}}
\newcommand{\Ti}[1]{T^{(#1)}}
\newcommand{\Fi}[1]{F^{(#1)}}
\newcommand{\ki}[1]{k^{(#1)}}

\newcommand{\yi}[1]{y^{(#1)}}
\newcommand{\zi}[1]{z^{(#1)}}
\newcommand{\si}[1]{s^{(#1)}}
\newcommand{\wi}[1]{w^{(#1)}}

\newcommand{\dtdp}[1]{\tau_{#1}}

\newcommand{\intgr}{\int_{t_{i}}^{t_{i+1}}}

\newcommand{\Ii}{I_{i}^{i+1}}
\newcommand{\Ji}{J_{i}^{i+1}}
\newcommand{\Ki}{K_{i}^{i+1}}

\begin{document}
	
\begin{titlepage}
	\begin{center}
		Simulation-Based Engineering Lab \\
		University of Wisconsin-Madison \\
		\LARGE Technical Report TR--2019--01\\
		\vfill 
		\LARGE Sensitivity Analysis \\ for Hybrid Systems and Systems with Memory \\
		\vspace{0.5in}
		\Large Radu Serban \\\medskip
		{\normalsize Department of Mechanical Engineering\\University of Wisconsin -- Madison} \\
		\vspace{0.2in}
		\Large Antonio Recuero \\\medskip
		{\normalsize Virtual Capability Department\\The Goodyear Tire \& Rubber Company} \\
		\vfill		
		\today		
	\end{center}
\end{titlepage}


\newpage
\begin{abstract}
We present an adjoint sensitivity method for hybrid discrete -- continuous systems, extending previously published forward sensitivity methods. We treat ordinary differential equations and differential-algebraic equations of index up to two (Hessenberg) and provide sufficient solvability conditions for consistent initialization and state transfer at mode switching points, for both the sensitivity and adjoint systems. Furthermore, we extend the analysis to so-called {\em hybrid systems with memory} where the dynamics of any given mode depend explicitly on the states at the last mode transition point. 
We present and discuss several numerical examples, including a computational mechanics problem based on the so-called exponential model constitutive material law for steel reinforcement under cyclic loading.

\vspace{0.25in}

\end{abstract}

\newpage
\tableofcontents
\newpage

\section{Introduction}\label{s:introduction}

Sensitivity analysis, the study of how variation in a model output can be apportioned to different sources of variation, has a wide range of applications in science and engineering.  These include model evaluation (finding the most and least influential parameters), generation of reduced-order models (reducing model complexity while preserving the input-output behavior), uncertainty quantification (characterizing and reducing uncertainty in model predictions), data assimilation (merging observations into a model to improve its accuracy), and dynamic optimization (finding model parameters for improved system response).

There is a significant body of work on methods and applications of sensitivity analysis, spanning at least the last three decades.  For dynamic models described by systems of ordinary differential equations (ODEs) or differential-algebraic equations (DAEs), like those considered herein, two main approaches for continuous 1\textsuperscript{st} order sensitivity analysis have been proposed and investigated.  
The forward sensitivity method (FSA)~\cite{caracotsios1985, maly1996numerical, feehery1997efficient}, which uses linearizations of the original model to calculate the state sensitivities, is relatively simple to formulate and implement, but can quickly become prohibitively expensive for large systems and large numbers of model parameters. 
Adjoint sensitivity analysis (ASA) provides an attractive alternative for problems with an arbitrary number of models parameters but relatively few output functionals for which gradients are sought.  More complicated, both from a formulation and implementation point of view, ASA requires the formulation of adjoint models~\cite{cacuci2003sensitivity} and, in the case of time-dependent problems, backward in time integration of a system that may depend on the forward solution. The adjoint models and their properties have been studied for systems of ODEs, index-1 and index-2 DAEs~\cite{cao2002adjoint, cao2003adjoint}, as well as for the index-3 DAEs in multibody dynamics~\cite{schaffer2005stability}.
%
%
General-purpose software for FSA and ASA of systems described by ODEs and DAEs is available and widely used~\cite{tolsma2000daepack, li2000software, hindmarsh2005sundials}.

While many physical phenomena can be described with continuous differential models, there are applications that require models for interaction between continuous and discrete phenomena. Such hybrid discrete -- continuous systems, with the continuous part described by systems of ODEs or DAEs and the discrete part modeled by finite automata, arise in control systems (e.g., safety interlock systems), biology (cell signaling models), chemical engineering (plant processing models), electrical engineering (heterogeneous electrical grid models), as well as mechanical engineering. In the area of computational mechanics, two applications stand out: collision and contact problems formulated in an event-driven framework; and rate-independent hysteretic phenomena in constitutive equations.
Forward sensitivity analysis for hybrid systems of this type was treated in a comprehensive manner by Barton and co-workers~\cite{galan1999parametric,tolsma2002hidden}, who highlighted the problems associated with hidden discontinuities during sensitivity analysis and the importance of properly computing and incorporating potential sensitivity jumps.
Adjoint methods for hybrid discrete -- continuous problems have seen relatively less attention. Recent work by Corner {\em et.al.}~\cite{sandu_2018_arXiv} treats ASA for hybrid multibody dynamical systems in the context of non-smooth contact dynamics. In~\cite{anitescuDiscreteASA}, the discrete adjoint method (adjoint of the discretized time-dependent problem) was used for computing derivative information for optimization of power system dynamics.

Here, we present an adjoint sensitivity method for hybrid discrete -- continuous systems, extending the forward sensitivity methods presented in~\cite{galan1999parametric, tolsma2002hidden}. We treat ordinary differential equations (ODEs) and DAEs of index up to two (Hessenberg) and provide sufficient solvability conditions for consistent initialization and state transfer at mode switching points, for both the sensitivity and adjoint systems. Furthermore, we extend the analysis to so-called {\em hybrid systems with memory} where the dynamics of any given mode depend explicitly on the states at the last mode transition point. 

This paper is structured as follows. In \S\ref{s:problem_formulation} we provide the framework for the type of hybrid systems considered herein. In \S\ref{s:SA_hybrid} we derive the sensitivity and adjoint systems for hybrid parameter-dependent problems whose dynamics are specified as index-0 (ODEs) or index-1 DAEs. The case of Hessenberg index-2 DAEs is treated separately in \S\ref{ss:HI2_hybrid}. In \S\ref{s:SA_memory} we discuss the sensitivity formulation for systems with memory, also providing an algorithm for propagating backward in time the adjoint systems while properly accounting for jumps at transition points. Numerical examples are discussed in \S\ref{s:examples}, including a computational mechanics problem based on the exponential model constitutive material law for steel structures under cyclic loading.

\section{Problem Formulation}\label{s:problem_formulation}

Let
\begin{equation}\label{e:DAE}
F \left( \dot x, x, p, t \right) = 0
\end{equation} 
be a parameter-dependent DAE system, where $x = \{ y, z\} \in \R^{N_x} = \R^{N_y + N_z}$, with $y$ and $\dot y$ being the differential states and their time derivatives, $z$ the algebraic variables, and $p \in \R^{N_p}$ a set of time-invariant problem parameters.

\begin{figure}
    \centering
    \includegraphics[width = 0.8\columnwidth]{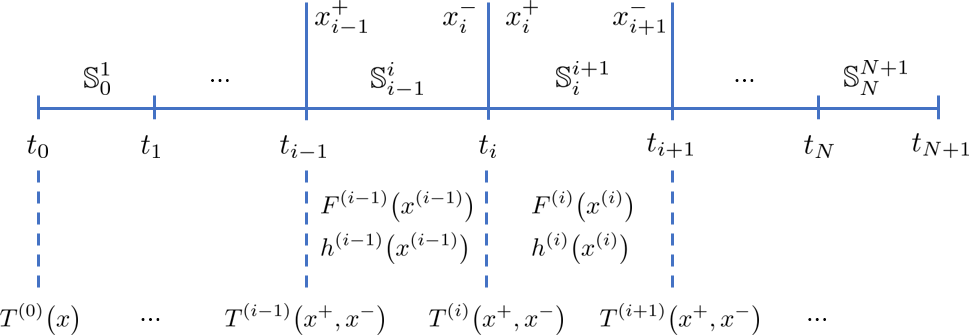}
    \caption{Formulation of hybrid discrete--continuous systems.}
    \label{f:hybrid}
\end{figure}

In this paper we consider hybrid discrete -- continuous systems which can be described following the formalism in~\cite{galan1999parametric,tolsma2002hidden}. Such systems can be in one of several modes $\Si$, $i=0,\ldots,N$ where each mode and the transitions between modes are characterized by (see also Fig.~\ref{f:hybrid}):
\begin{description}
\item{Variables.}
The state of the system in mode $\Si$ is described by a set of variables
$\left\{ \dot y^{(i)}, y^{(i)}, z^{(i)}\right\}$, functions of time $t$ and of the time-invariant parameters $p$;
\item{Dynamics.}
In mode $\Si$, the dynamics of the system are given by DAEs of the form of Eq.~\ref{e:DAE}. In this paper, we consider problems for which $\text{rank}\left([F_{\dot y} | F_z]\right) = N_y + N_z$ which is true for index-0 (ODEs) and most index-1 DAEs.  The case of different DAE structures (namely Hessenberg index-2) is treated separately in \S\ref{ss:HI2_hybrid}.
\item{Transitions.}
The system transitions to a different mode based on the so-called {\em transition conditions} $h^{(i-1)}\left(\dot y^{(i-1)}, y^{(i-1)}, z^{(i-1)}, p, t\right)$. These conditions implicitly define the transition time $t_{i}, i = 1,\ldots,N$ at which the system switches to mode $\Si$.  

The {\em transition functions} $\Ti{i}\left(\dot y^{(i)}, y^{(i)}, z^{(i)} , \dot y^{(i-1)}, y^{(i-1)}, z^{(i-1)} , p,t\right)$ define a system of equations that map the final states in mode $\mathbb{S}_{i-1}^{i}$  (at the transition time $t_i$) to the initial states in the next mode $\Si$ (after the transition time). 
\end{description}
Without loss of generality, we assume that the initial and final simulation times do not depend on the problem parameters $p$. Initial conditions are provided as:
\begin{equation}
    \Ti{0}(\dot y, y, p, t) = 0 \,,
\end{equation}
which are assumed to provide $N_y$ additional conditions consistent with the DAE of Eq.~\ref{e:DAE}. In other words,
\begin{equation*}
    \begin{bmatrix}
    \Ti{0}_{\dot y} & \Ti{0}_y & \Ti{0}_z \\
    \Fi{0}_{\dot y} & \Fi{0}_y & \Fi{0}_z 
    \end{bmatrix} \text{ is assumed to be nonsingular at } t=t_0 \,,
\end{equation*}
where $\Fi{0}$ defines the system dynamics on the first mode $\mathbb{S}_0^1$.

With the above description, we assume that the system undergoes $N$ transitions at times $t_i, i=1,\ldots,N$ and evolves through a sequence of $N+1$ modes $\Si$, each over $t \in [t_i, t_{i+1}]$. 
At each transition time $t_i$, we denote by $a^-_i$ a quantity $a$ evaluated {\em before} the transition; e.g., $y^-_i = y(t^-_i) = y(t_i - 0)$ and by $a^+_i$ the same quantity {\em after} the transition (possibly including a discontinuity); e.g, $y^+_i = y(t^+_i) = y(t_i + 0)$.

The transition times are defined implicitly by zero-crossings of the transition conditions $\hi{i-1}$ as solution of the equations $\hi{i-1}\left(\dot y^-_i , y^-_i, z^-_i, p, t_i\right) = 0$ and the initial values for the DAE problem over the next mode are dictated by the system of equations $\Ti{i}\left(\dot y^+_i , y^+_i, z^+_i, \dot y^-_i , y^-_i, z^-_i, p, t_i\right) = 0$.

Solvability conditions of the resulting hybrid problem at each transition point stem from imposing, at each transition time $t_i$, the following nonlinear equations:
\begin{subequations}\label{e:transition_hybrid}
\begin{align}
    &\dot y^-_i = \dot y^{(i-1)} (p, t_i) \,,\quad y^-_i = \yi{i-1}(p, t_i) \,,\quad z^-_i = \zi{i-1}(p, t_i) \,, \\
    &\hi{i-1}\left(\dot y^-_i , y^-_i, z^-_i, p, t_i\right) = 0 \,, \\
    &\Ti{i}\left(\dot y^+_i , y^+_i, z^+_i, \dot y^-_i , y^-_i, z^-_i, p, t_i\right) = 0 \,, \\
    &\Fi{i}\left(\dot y^+_i , y^+_i, z^+_i, p, t_i\right) = 0 \,,
\end{align}
\end{subequations}
where $\Fi{i}(\dot y, y, z, p, t)$ defines the system dynamics over the next mode, $\Si$.
Sufficient conditions for the solvability of this problem (based on the implicit function theorem) are
\begin{equation}\label{e:solvability_hybrid}
\begin{split}
&\hi{i-1}_{\dot y} \ddot y^-_i +
\hi{i-1}_{y} \dot y^-_i +
\hi{i-1}_{z} \dot z^-_i +
\hi{i-1}_{t} \ne 0 \\
\text{and}\qquad& \\
&\begin{bmatrix}
\Ti{i}_{\dot y^+} & \Ti{i}_{y^+} & \Ti{i}_{z^+} \\
\Fi{i}_{\dot y} & \Fi{i}_{y} & \Fi{i}_{z}
\end{bmatrix} \text{ nonsingular at } t=t_i \,,
\end{split}
\end{equation}
with the second simply being the conditions for consistent DAE initialization at $t_i$. Note that all partial derivatives above must be evaluated at $t_i$.


\section{Sensitivity Analysis for Hybrid Systems}\label{s:SA_hybrid}

Forward sensitivity analysis and the general equations for parameteric sensitivity functions for hybrid systems like those described in \S\ref{s:problem_formulation} have been presented in ~\cite{galan1999parametric, tolsma2002hidden}.  We briefly review these results in \S\ref{ss:FSA_hybrid} insofar as they are a required step in deriving the adjoint equations and performing adjoint sensitivity analysis for such systems, as presented in \S\ref{ss:ASA_hybrid}.


\subsection{Sensitivity Equations (FSA)}\label{ss:FSA_hybrid}
Differentiating the system DAE of Eq.~\ref{e:DAE} with respect to the independent variables $p$, the sensitivity trajectories in a given mode $\Si$ can be obtained by integrating the linear time varying sensitivity DAE simultaneously with the original DAE:
\begin{equation}\label{e:FSA}
    F_{\dot y} \dot s + F_y s + F_z w + F_p = 0 \,,
\end{equation}
where $s$ represents the $N_y \times N_p$ matrix of differential state sensitivities and $w$ is the $N_z \times N_p$ matrix of algebraic state sensitivities.

At each transition point $t_i$, in addition to computing consistent initial values for the DAE states in the new mode, 
the sensitivity system in the new mode must also be consistently initialized -- often requiring calculating explicit jumps in the sensitivities.  This procedure requires first calculating the sensitivity of the transition times themselves which can be obtained from a formal differentiation of the transition conditions of Eq.~\ref{e:transition_hybrid} to obtain:
\begin{subequations}\label{e:transition_hybrid_FSA}
\begin{align}
&\dot s^-_i = \dot s^{(i-1)}(p, t_i) \,,\quad s^-_i = \si{i-1}(p, t_i) \,,\quad w^-_i = \wi{i-1}(p, t_i) \,,
\\
&\hi{i-1}_{\dot y} \left(\dot s^-_i + \ddot y^-_i \dtdp{i} \right)+
\hi{i-1}_{y} \left(s^-_i + \dot y^-_i \dtdp{i}\right)+
\hi{i-1}_{z} \left(w^-_i + \dot z^-_i \dtdp{i}\right)+
\hi{i-1}_{t} \dtdp{i} +
\hi{i-1}_{p} = 0 \,,
\\
\begin{split}
&\Ti{i}_{\dot y^+} \left(\dot s^+_i + \ddot y^+_i \dtdp{i} \right)+
\Ti{i}_{y^+} \left(s^+_i + \dot y^+_i \dtdp{i}\right)+
\Ti{i}_{z^+} \left(w^+_i + \dot z^+_i \dtdp{i}\right) \\
&\qquad+\Ti{i}_{\dot y^-} \left(\dot s^-_i + \ddot y^-_i \dtdp{i} \right)+
\Ti{i}_{y^-} \left(s^-_i + \dot y^-_i \dtdp{i}\right)+
\Ti{i}_{z^-} \left(w^-_i + \dot z^-_i \dtdp{i}\right) \\
&\qquad+\Ti{i}_{t} \dtdp{i} + \Ti{i}_{p} = 0 \,,
\end{split}
\\
&\Fi{i}_{\dot y} \left(\dot s^+_i + \ddot y^+_i \dtdp{i}\right) +
\Fi{i}_y \left(s^+_i + \dot y^+_i \dtdp{i}\right) +
\Fi{i}_z \left(w^+_i + \dot z^+_i \dtdp{i}\right) +
\Fi{i}_t \dtdp{i} +
\Fi{i}_p = 0 \,,
\end{align}
\end{subequations}
where $\tau_i = {dt_i}/{dp}$ is the sensitivity of the transition time $t_i$ with respect to the problem parameters $p$, $s^{(i-1)}(p, t)$ and $w^{(i-1)}(p, t)$ are the state sensitivity trajectories over the previous mode, and all partial derivatives are evaluated at $t_i$.  It can be easily seen that sufficient conditions for the solvability of the above system of equations are the same as those of Eq.~\ref{e:solvability_hybrid}.

Note that $\ddot y$ and $\dot z$ in Eq.~\ref{e:transition_hybrid_FSA} can be obtained by differentiating Eq.~\ref{e:DAE} with respect to the independent variable $t$:
$F_{\dot y} \ddot y + F_y \dot y + F_z \dot z + F_t = 0$, which is solvable under the assumption that $[F_{\dot y} | F_z]$ has full rank.

With the above, we can propagate forward in time the augmented DAE and sensitivity system through the sequence of modes $\Si, i = 0,\ldots,N$, properly taking into account any required discrete updates (jumps) at the transition points $t_i, i=1,N$.


\subsection{Adjoint Sensitivity Analysis (ASA)}\label{ss:ASA_hybrid}

The augmented DAE and sensitvity system of Eqs.~\ref{e:DAE} and \ref{e:FSA} can be quite large, especially for large problems.  Even though there are efficient algorithms that exploit the structure of this augmented system~\cite{maly1996numerical,feehery1997efficient,hindmarsh2005sundials}, the cost of such computations can become prohibitive when there are many problem parameters.  The adjoint sensitivity method has proven to be a viable, effective alternative, at least when interested in the parameter sensitivity of one or only a few functionals that depend on the system state~\cite{cacuci2003sensitivity,cao2002adjoint,cao2003adjoint,hindmarsh2005sundials}.

In this section, we present the derivation of the adjoints for hybrid systems and provide the necessary transition conditions for the adjoint variables, required for correct ASA of such problems. As before, we consider here index-0 (i.e., ODEs) and index-1 problems; the case of the Hessenberg index-2 DAE structure is treated in \S\ref{ss:HI2_hybrid}.

Consider a functional of the system state trajectory of the form:
\begin{equation}\label{e:G_hybrid}
    G(p) = \int_{t_0}^{t_f} g(y,z,p,t) dt = \sum_{i=0}^{N} \intgr g(y,z,p,t) dt \,.
\end{equation}
Its gradient with respect to the time-invariant problem parameters is then
\begin{equation}
    \frac{dG}{dp} = \sum_{i=0}^{N} \left[
    \intgr\left(g_y s + g_z w + g_p\right) dt + 
    g\left(y^-_{i+1},z^-_{i+1},p,t_{i+1}\right) \dtdp{i+1} -
    g\left(y^+_i,z^+_i,p,t_i\right)\dtdp{i}
    \right] \,,
\end{equation} 
where we have used the Leibniz integral rule to account for the dependency of the transition times on $p$ (recall that $\dtdp{i} = {dt_i}/{dp}$).
Introducing adjoint variables $\lambda \in \R^{N_y + N_z}$ and taking into consideration Eq.~\ref{e:FSA},
\begin{equation}
\begin{split}
\frac{dG}{dp} = \sum_{i=0}^{N} &\left[
\intgr\left(g_y s + g_z w + g_p\right) dt + g^-_{i+1} \dtdp{i+1} - g^+_i \dtdp{i} \right. \\
&+ \left.\intgr \lambda^T \left(F_{\dot y} \dot s + F_y s + F_z w + F_p \right) dt\right] \,.
\end{split}
\end{equation}
Using integration by parts
\begin{equation}
\intgr \lambda^T F_{\dot y} \dot s dt
= \left.\left(\lambda^T F_{\dot y} s\right)\right|_{t_i}^{t_{i+1}}
- \intgr \frac{d}{dt}\left(\lambda^T F_{\dot y}\right) s dt \,,
\end{equation}
and rearranging terms, we obtain
\begin{equation}\label{e:dGdp1}
\begin{split}
\frac{dG}{dp} = \sum_{i=0}^{N} &\left[
\intgr \left( g_p + \lambda^T F_p\right) dt +
g^-_{i+1}\dtdp{i+1} - g^+_i \dtdp{i} \right. \\
& + \left.\left(\lambda^-_{i+1}\right)^T F_{\dot y}\left(t^-_{i+1}\right) s^-_{i+1} -
\left(\lambda^+_{i}\right)^T F_{\dot y}\left(t^+_{i}\right) s^+_{i}
\right] \,,
\end{split}
\end{equation}
if the adjoint variables are selected to satisfy
\begin{equation}\label{e:ASA}
\begin{split}
&-\frac{d}{dt}\left(\lambda^T F_{\dot y}\right) + \lambda^T F_y + g_y = 0 \\
&\lambda^T F_z + g_z = 0 \,.
\end{split}    
\end{equation}

The computation of the gradient ${dG}/{dp}$ using the adjoint method is then completed if we impose appropriate jump conditions on $\lambda$ at the transition times $t_i$ so that the sensitivities $s_i$ cancel out in Eq.~\ref{e:dGdp1}.  
To obtain these conditions, we start by rearranging terms and favorably changing the summation indices:
\begin{equation}\label{e:dGdp2}
\begin{split}
\frac{dG}{dp}& = 
\sum_{i=0}^{N} \left[\intgr \left( g_p + \lambda^T F_p\right) dt \right] \\
&\quad 
+ \sum_{i=0}^N \left[ g^-_{i+1}\dtdp{i+1} - g^+_i \dtdp{i}
+ \left(\lambda^-_{i+1}\right)^T F_{\dot y}\left(t^-_{i+1}\right) s^-_{i+1}
- \left(\lambda^+_{i}\right)^T F_{\dot y}\left(t^+_{i}\right) s^+_{i}
\right] \\
& = \sum_{i=0}^{N} \left[\intgr \left( g_p + \lambda^T F_p\right) dt \right] \\
&\quad
+ \lambda_f^T F_{\dot y}\left(t_f\right) s_f
- \lambda_0^T F_{\dot y}\left(t_0\right) s_0 \\
&\quad
+ \sum_{i=1}^N \left[
\left( g^-_i - g^+_i \right) \dtdp{i} 
+ \left(\lambda^-_{i}\right)^T F_{\dot y}\left(t^-_{i}\right) s^-_{i}
- \left(\lambda^+_{i}\right)^T F_{\dot y}\left(t^+_{i}\right) s^+_{i}
\right] \,,
\end{split}  
\end{equation}
where we used the fact that the initial ($t_0$) and final ($t_f = t_{N+1}$) times do not depend on the problem parameters $p$. For index-0 and index-1 DAEs, the initial conditions (at $t = t_f$) for the adjoint variables can be simply 
\begin{equation}
    \lambda_f^T F_{\dot y} = 0 \,.
\end{equation}
As shown in \S\ref{ss:HI2_hybrid} this choice will not suffice for a Hessenberg index-2 system.

Finally, the jump conditions for the adjoint variables at each transition time are inferred from the components of the last term in Eq.~\ref{e:dGdp2} by imposing
\begin{equation*}
    \left[ 
    \left(\lambda^-_{i}\right)^T F_{\dot y} +
    \left( g^-_i - g^+_i \right) \beta_i
    \right] s^-_{i} =
    \left(\lambda^+_{i}\right)^T F_{\dot y}\left(t^+_{i}\right) s^+_{i} \,,
\end{equation*}
where we formally expressed $\dtdp{i} = \alpha_i + \beta_i s^-_i$ from Eq.~\ref{e:transition_hybrid_FSA}. These relations, together with the algebraic constraint in Eq.~\ref{e:ASA}, allow us to calculate $\lambda^-_i$, the adjoint variables after crossing a transition point, in terms of $\lambda^+_i$, the adjoint variables at the end of (backward) integration over the previous mode. Note that the assumption $[F_{\dot y} | F_z]$ nonsingular again provides sufficient solvability conditions.


\subsection{Hessenberg index-2 problems}\label{ss:HI2_hybrid}

The assumptions used previously for establishing the solvability conditions at transition points for the original DAE (Eq.~\ref{e:DAE}), the linear sensitivity system (Eq.~\ref{e:FSA}), and the adjoint problem (Eq.~\ref{e:ASA}) will not suffice for a Hessenberg index-2 DAE system.  We provide these derivations in this section. The main reason for including this type of problems in our analysis is that this is the form of DAEs obtained when modeling incompressible flow by discretized Navier-Stokes equations in CFD~\cite{ascher1998computer} or with a stabilized index reduction of the index-3 equations of motion in multibody system problems~\cite{gear1985automatic, brenan1996numerical}.

Let the dynamics of a hybrid system such as those defined in \S\ref{s:problem_formulation} be described, in each mode $\Si$, by DAEs in {\em Hessenberg index-2} form:
\begin{subequations}\label{e:Hi2}
\begin{align}
\dot y &= f^{(i)}(y,z,p,t) \label{e:Hi2_diff} \\
0 &= k^{(i)}(y,p,t) \,, \label{e:Hi2_alg}
\end{align}
\end{subequations}
where $f^{(i)}$ is the right-hand side of the differential equations, $k^{(i)}$ denotes the algebraic equations, and $C^{(i)}B^{(i)}$ is nonsingular, with $A^{(i)} = f^{(i)}_y$, $B^{(i)}=f^{(i)}_z$, and $C^{(i)} = k^{(i)}_y$.

Obtaining consistent initial conditions for higher-index DAEs is more complicated and no general recipe exists.  It is clear that not any choice of $N_y$ initial conditions for the differential variables will do, as $y_0 = y(t_0)$ must satisfy the algebraic equations~\ref{e:Hi2_alg}.  Furthermore, higher-index DAEs have so-called {\em hidden manifolds} to which the solution is restricted; for the system above, additional constraints are given by the differentiated algebraic equations: $C f + k_t = 0$.  See~\cite{EstevezSchwarz2009} for further details on consistent initialization of DAE and section \S 4.2 in~\cite{li2000software} for a discussion of practical implementation considerations.
Here, we assume that initial conditions $y_0$ consistent with Eq.~\ref{e:Hi2} are provided.  Furthermore, we assume these are given implicitly as $0 = \Ti{0}(y, p, t) \in \R^{N_y-N_z}$.  A sufficient condition (based on the implicit function theorem) for the solvability of the resulting nonlinear system is then
\begin{equation*}
    \begin{bmatrix}\Ti{0}_y\\C^{(0)}\end{bmatrix} \text{ nonsingular at } t = t_0\,.
\end{equation*}
Initial conditions for the algebraic variables are then obtained by solving the nonlinear system $C^{(0)} f^{(0)} + k^{(0)}_t = 0$, evaluated at $\{t_0, y_0\}$.  Note that the regularity condition ($CB$ nonsingular) is critical in providing solvability conditions for this step.  Finally, set $\dot y_0 = f^{(0)}(y_0, z_0, p, t_0)$.

As before, consider a hybrid discrete -- continuous system whose dynamics in mode $\Si$ are given by equations of the form~\ref{e:Hi2}, with transitions\footnote{For simplicity, we ignore here more complex transition conditions and functions which, in general, may also involve the time derivatives and/or algebraic variables. The case treated here is typical for problems in CFD or multibody dynamics.} defined by the transition conditions $\hi{i}\left(\yi{i}, \zi{i}, p, t\right)$ and the associated transition functions $\Ti{i}\left(\yi{i}, \yi{i-1}, p, t\right)$.

The state transition conditions of the resulting hybrid problem become
\begin{subequations}\label{e:transition_Hi2_hybrid}
\begin{align}
    &y^-_i = \yi{i-1}(p, t_i) \,,\quad z^-_i = \zi{i-1}(p, t_i) \,, \\
    &\hi{i-1}\left(y^-_i, z^-_i, p, t_i\right) = 0 \,, \label{e:transition_Hi2_time}\\
    &\Ti{i}\left(y^+_i, y^-_i, p, t_i\right) = 0 \,, \label{e:transition_Hi2_state}\\
    &\dot y^+_i = f^{(i)}\left( y^+_i, z^+_i, p, t_i \right) \,, \label{e:transition_Hi2_diff}\\
    &0 = k^{(i)}(y^+_i, p, t_i)\,, \label{e:transition_Hi2_alg}\\
    &0 = C^{(i)} \dot y^+_i + k^{(i)}_t \,. \label{e:transition_Hi2_hidden}
\end{align}
\end{subequations}
Sufficient conditions for computing the state transitions are:
\begin{equation}\label{e:solvability_Hi2_hybrid}
\begin{split}
&\hi{i-1}_{y} \dot y^-_i +
\hi{i-1}_{z} \dot z^-_i +
\hi{i-1}_{t} \ne 0 \\
\text{and}\qquad& \\
&\begin{bmatrix}\Ti{i}_{y^+}\\C^{(i)}\end{bmatrix} \text{ nonsingular at } t=t_i \,,
\end{split}
\end{equation}
in addition to the regularity condition on $\Si$, namely $\exists \left(C^{(i)}B^{(i)}\right)^{-1}$.  The transition time $t_i$ is solved from Eq.~\ref{e:transition_Hi2_time}; the differential states $y^+_i$ are obtained by simultaneously solving Eqs.~\ref{e:transition_Hi2_state} and \ref{e:transition_Hi2_alg}; finally, the algebraic states after the transition are obtained from the nonlinear system of Eq.~\ref{e:transition_Hi2_hidden} and the state derivatives $\dot y^+_i$ calculated with Eq.~\ref{e:transition_Hi2_diff}.

The corresponding sensitivity equations are
\footnote{For clarity, we drop all subscripts and superscripts in the remainder of this section, except where necessary.}
\begin{subequations}\label{e:Hi2_FSA}
\begin{align}
    \dot s &= f_y s + f_z w + f_p \label{e:Hi2_FSA_diff} \\
    0 &= k_y s + k_p \label{e:Hi2_FSA_alg}
\end{align}
\end{subequations}
The sensitivity transition conditions are obtained in a manner completely analogous to Eq.~\ref{e:transition_hybrid_FSA} and are not reproduced here. Similar to the case of index-0 and index-1 systems, the solvability conditions for computing the sensitivity of the transition times and the state sensitivity jump transitions are those of Eq.~\ref{e:solvability_Hi2_hybrid}.
Furthermore, we can formally write $\tau_i = \alpha_i + \beta_i s^-_i$ and $s^+_i = \Gamma_i + \Delta_i s^-_i$.

Considering the functional $G$ of Eq.~\ref{e:G_hybrid} and following a derivation similar to that of \S\ref{ss:ASA_hybrid}, its gradient with respect to the problem parameters $p$ is:
\begin{equation}\label{e:dGdp2_Hi2}
\begin{split}
\frac{dG}{dp}& = 
\sum_{i=0}^{N} \left[\intgr \left( g_p + \lambda^T f_p + \mu^T k_p\right) dt \right] \\
&\quad - \lambda_f^T s_f + \lambda_0^T s_0 \\
&\quad
+ \sum_{i=1}^N \left[
\left( g^-_i - g^+_i \right) \dtdp{i} - \left(\lambda^-_{i}\right)^T s^-_{i} + \left(\lambda^+_{i}\right)^T s^+_{i}
\right] \,,
\end{split}  
\end{equation}
where the adjoint variables $\lambda$ and $\mu$ were chosen to satisfy
\begin{subequations}\label{e:Hi2_Adjoint}
\begin{align}
\dot \lambda^T + \lambda^T f_y + \mu^T k_y + g_y &= 0 \label{e:Hi2_Adjoint_diff} \\
\lambda^T f_z + g_z &= 0 \label{e:Hi2_Adjoint_alg} \,,
\end{align}
\end{subequations}
in each interval $[t_i, t_{i+1}]$. Because of the structure of the DAE, the choice $\lambda_f = 0$ as initial conditions for the adjoint hybrid system will not suffice, as this could conflict with Eq.~\ref{e:Hi2_Adjoint_alg} if $g(y,z,p,t)$ depends explicitly on $z$.  Instead, with all quantities evaluated at $t = t_f$ using functions corresponding to mode $\mathbb{S}_N^{N+1}$, we require that
\begin{equation*}
    \lambda_f^T = \xi^T C \,,
\end{equation*}
for some $\xi$ yet to be determined~\cite{cao2003adjoint}. Inserting into Eq.~\ref{e:Hi2_Adjoint_alg}, we have
\begin{equation*}
\xi^T C B = -g_z \quad\Rightarrow\quad \xi^T = -g_z \left( C B\right)^{-1} \,.
\end{equation*}
Therefore,
\begin{equation*}
\lambda_f^T = -g_z \left( C B\right)^{-1} C
\end{equation*}
and since $C s_f = -k_p $ (from Eq.~\ref{e:Hi2_FSA_alg})
\begin{equation*}
\lambda^T_f s_f = g_z \left( C B\right)^{-1} k_p \,.
\end{equation*}
Consistent initial conditions for $\mu_f$ are obtained from the hidden constraint resulting from differentiation of the algebraic equations~\ref{e:Hi2_Adjoint_alg}. Indeed,
\begin{equation*}
    \lambda^T B + g_z = 0 \quad\Rightarrow\quad \dot\lambda^T B + \lambda^T \dot B + \dot g_z = 0
\end{equation*}
and inserting $\dot\lambda$ from Eq.~\ref{e:Hi2_Adjoint_diff}
\begin{equation*}
    \mu_f^T = \left(\lambda_f ^T (\dot B - AB) + ( \dot g_z - g_y B ) \right) \left(CB\right)^{-1} \,.
\end{equation*}

A similar procedure can be applied to obtain the transfer conditions (jumps) of the adjoint variables at each transition point.
Indeed, we are seeking $\lambda^-_i$ such that each term in the last sum on the right-hand side of Eq.~\ref{e:dGdp2_Hi2} can be evaluated without the need for computing the state sensitivities $s^-_i$ and $s^+_i$.  Using $\tau_i = \alpha_i + \beta_i s^-_i$ and $s^+_i = \Gamma_i + \Delta_i s^-_i$, after simple algebraic manipulations, we have 
\begin{equation*}
\left( g^-_i - g^+_i \right) \dtdp{i} - \left(\lambda^-_{i}\right)^T s^-_{i} + \left(\lambda^+_{i}\right)^T s^+_{i}
= \left( g^-_i - g^+_i \right)\alpha_i + \left(\lambda^+_{i}\right)^T \Gamma_i
- \left[ \left(\lambda^-_{i}\right)^T - X_i^T \right] s^-_i \,,
\end{equation*}
where
\begin{equation*}
X_i^T = \left( g^-_i - g^+_i \right) \beta_i + \left(\lambda^+_{i}\right)^T \Delta_i \,.
\end{equation*}
Note that we cannot simply impose $\lambda^-_i = X_i$ as this could conflict with the algebraic equations in Eq.~\ref{e:Hi2_Adjoint_alg}.  Instead, using quantities evaluated at $t^-_i$, we seek a vector $\zeta_i$ such that
\begin{equation*}
\left( \lambda^-_i - X_i\right)^T = \zeta_i^T C \,.
\end{equation*}
Inserting into Eq.~\ref{e:Hi2_Adjoint_alg} (evaluated at $t^-_i$ using the DAE for mode $\mathbb{S}_{i-1}^i$), we obtain
\begin{equation*}
\left( \zeta_i^T C + X_i^T \right) B = - g_z
\quad\Rightarrow\quad
\zeta_i^T = -\left(g_z + X_i^T B \right) \left( C B\right)^{-1} \,.
\end{equation*}
Therefore,
\begin{equation*}
\left(\lambda^-_i\right)^T = X_i^T - \left(g_z + X_i^T B \right) \left( C B\right)^{-1} C
\end{equation*}
and since $C s^-_i = -k_p $ (from Eq.~\ref{e:Hi2_FSA_alg})
\begin{equation*}
\left[ \left(\lambda^-_{i}\right)^T - X_i^T \right] s^-_i =
\left(g_z + X_i^T B \right) \left( C B\right)^{-1} k_p \,.
\end{equation*}
The gradient of Eq.~\ref{e:dGdp2_Hi2} then becomes
\begin{equation}
\begin{split}
\frac{dG}{dp}& = 
\sum_{i=0}^{N} \left[\intgr \left( g_p + \lambda^T f_p + \mu^T k_p\right) dt \right] \\
&\quad 
+ \lambda_0^T s_0 - \left. g_z \left( C B\right)^{-1} k_p \right|_{t=t_f}  \\
&\quad
+ \sum_{i=1}^N \left[
\left( g^-_i - g^+_i \right)\alpha_i + \left(\lambda^+_{i}\right)^T \Gamma_i - 
\left. \left(g_z + X_i^T B \right) \left( C B\right)^{-1} k_p \right|_{t=t^-_i}
\right] \,.
\end{split}
\end{equation}


\section{Sensitivity Analysis for Hybrid Systems with Memory}\label{s:SA_memory}

In some situations, the system dynamics in mode $\Si$ depends explicitly on the states at the transition into mode $\Si$.  Such a problem is the computational mechanics example discussed in \S\ref{ss:example_EM_memory} which uses the so-called {\em exponential model}~\cite{vaiana2018class} model for steel reinforcement bars under cyclic loading. 

In such cases, this additional explicit dependency must be properly taken into account in deriving the sensitivity and adjoint systems and the corresponding transition rules. For brevity in exposition, we consider here only hybrid continuous-discrete systems whose dynamics between any two transition points are described by DAEs in {\em Hessenberg index-1} form:
\begin{subequations}\label{e:DAE_mem}
\begin{align}
\dot y & = f^{(i)}(y, z, y^* , z^*, p) \label{e:DAE_mem_diff}\\
0 & = k^{(i)}(y, z, y^*, z^*, p) \label{e:DAE_mem_alg} \,,
\end{align}
\end{subequations}
over the interval $[t_0, t_f]$, with $k_z$ nonsingular and $p \in \R^{N_p}$ a set of time-independent problem parameters.  The differential variables are $y \in \R^{N_y}$, $z \in \R^{N_z}$ are the algebraic variables, and without loss of generality we assume that $f$ and $k$ do not explicitly depend on time. Furthermore, we assume that consistent initial conditions are provided and, for simplicity, we assume that $y_0 \equiv y(t_0)$ does not depend on $p$.
Let the transition times $t_i$ between consecutive modes be implicitly defined by the transition conditions\footnote{Without loss of generality, we assume that the transition conditions $h^{(i-1)}$ do not depend explicitly on the algebraic variables, as those can always be formally expressed as $z = z(y)$ from the algebraic equations~\ref{e:DAE_mem_alg}.}
\begin{equation}\label{e:trans_cond_mem}
\hi{i-1}(y^-_i, p) = 0 \,, \quad i = 1,\ldots,N \,.
\end{equation}
Let $t_f = t_{N+1}$ and let $\Si$ be the mode over the interval $[t_{i}, t_{i+1}]$.  The transition into $\Si$ at $t_i$ is defined by the transition functions
\begin{equation}\label{e:trans_fct_mem}
\Ti{i}(y^+_i, z^+_i, y^-_i, z^-_i, p) = 0 \,,
\end{equation}
which, in general, may also include explicit dependencies on time or the time derivatives ${\dot y}^-_i$ and ${\dot y}^+_i$, before and after the transition, respectively. Equation~\ref{e:trans_fct_mem} implicitly defines $y^+_i$ and $z^+_i$, the system states after the transition. 
While in mode $\Si$, the equations of motion are those of Eq.~\ref{e:DAE_mem} with $y^* = y^+_{i} \equiv y(t^+_{i})$ and $z^* = z^+_{i} \equiv z(t^+_{i})$.
%
Solvability of the hybrid problem requires that the transition conditions of Eq.~\ref{e:trans_fct_mem} be consistent with the algebraic equations of Eq.~\ref{e:DAE_mem_alg} at time $t^+_i$. A sufficient condition is:
\begin{equation*}
\begin{bmatrix}
\Ti{i}_{y^+} & \Ti{i}_{z^+} \\ \ki{i}_{y}(t^+_i) & \ki{i}_{z}(t^+_i)
\end{bmatrix} \text{ full rank. }
\end{equation*}
In many physical phenomena, including the example problem described in \S\ref{ss:example_EM_memory}, it is common for the system states to be $\mathcal{C}^0$ continuous at transition points; i.e., $y^+_i = y^-_i$ and $z^+_i = z^-_i$.

We call $y^*$ and $z^*$ {\em memory states} and the problem of Eq.~\ref{e:DAE_mem} a {\em system with memory}. In this section, we extend the previous formulation to perform (forward and adjoint) sensitivity analysis of such systems. The new elements in this derivation are related to the implicit dependency of the memory states on problem parameters.


\subsection{Sensitivity Equations (FSA)}\label{ss:FSA_memory}

Denoting by $s \equiv {dy}/{dp}$ and $w \equiv {dz}/{dp}$ the state sensitivities, a formal differentiation of Eq.~\ref{e:DAE_mem} leads to the following forward sensitivity equations:
\begin{subequations}\label{e:FSA_mem}
\begin{align}
\dot s & = f^{(i)}_y s + f^{(i)}_z w + f^{(i)}_{y^*} s^* + f^{(i)}_{z^*} w^* + f^{(i)}_p \label{e:FSA_mem_diff} \\
0 & = k^{(i)}_y s + k^{(i)}_z w + k^{(i)}_{y^*} s^* + k^{(i)}_{z^*} w^* + k^{(i)}_p \label{e:FSA_mem_alg} \,,
\end{align}	
\end{subequations}
where, in mode $\Si$, we have $s^* = s^+_{i} \equiv s(t^+_{i})$ and $w^* = w^+_{i} \equiv w(t^+_{i})$.

Consistent initialization of the FSA problem after a transition can be carried out as follows. Consider the transition at $t_i$ and, for clarity, drop the superscripts $(i)$ and subscripts $i$.
First, formal differentiation of the transition condition implies
\begin{equation}
h_y \left( s^- + \dot y^- \tau \right) + h_p = 0 \,,
\end{equation}
where we include the implicit dependency on $p$ through the time derivative $\dot y$. Recall that $\tau = dt_i/dp$.
This determines the sensitivity of the transition time as:
\begin{equation}\label{e:FSA_trans_time_mem}
\tau = - \frac{1}{h_y \dot y^-} \left( h_p + h_y s^- \right) \in \R^{1 \times N_p} \,,
\end{equation}
where all quantities are evaluated at $t^-_i$.

Next, a formal differentiation of the transition functions of Eq.~\ref{e:trans_fct_mem} gives
\begin{equation}\label{e:FSA_trans_fct_mem}
T_{y^+} \left( s^+ + \dot y^+ \tau \right) +
T_{z^+} \left( w^+ + \dot z^+ \tau \right) +
T_{y^-} \left( s^- + \dot y^- \tau \right) +
T_{z^-} \left( w^- + \dot z^- \tau \right) +
T_p = 0 \,,
\end{equation}
which, together with the consistency conditions of Eq.~\ref{e:FSA_mem_alg}, determines the sensitivity states after the transition, $s^+$ and $w^+$, given $s^-$ and $w^-$.
Note that the time derivatives of the algebraic states, $\dot z$, can be obtained simultaneously with $\ddot y$ by taking the time derivative of Eq.~\ref{e:DAE_mem} and solving, at time $t^+_i$, the resulting linear system with a nonsingular matrix
\begin{equation*}
\begin{bmatrix} I & -f_z \\ 0 & k_z\end{bmatrix}
\end{equation*}
(based on the Hessenberg index-1 condition).
In other words, $s^+$, $w^+$, and $\dot s^+$ are solved for simultaneously considering Eqs.~\ref{e:FSA_trans_fct_mem} and \ref{e:FSA_mem}; this problem is well posed given the solvability conditions for the transition functions of the original problem.

Integration of the forward sensitivity equations can then be carried out simultaneously with the original problem of Eq.~\ref{e:DAE_mem}, including the necessary jumps in sensitivity variables at each mode transition time.
 

\subsection{Adjoint Sensitivity Analysis (ASA)}\label{ss:ASA_memory}
 
Consider next a functional of the form
\begin{equation}\label{e:G_mem}
G(p) = \int_{t_0}^{t_f} g(y,z,p,t) dt = \sum_{i=0}^N \int_{t_i}^{t_{i+1}} g(y,z,p,t) dt
\end{equation}
and the problem of calculating the gradient of $G$ with respect to the problem parameters $p$.

Differentiating Eq.~\ref{e:G_mem} and using the Leibniz integral rule to take into account the dependency of the transition times $t_i$ on the problem parameters $p$, we have
\begin{equation}
\frac{dG}{dp} = \sum_{i=0}^N \left[ 
\intgr \left( g_y s + g_z w + g_p \right) dt + g^-_{i+1} \dtdp{i+1} - g^+_{i} \dtdp{i}
\right]
\end{equation}
where $g^-_{i+1} \equiv g(y^-_{i+1}, z^-_{i+1}, p, t_{i+1})$ and $g^+_{i} \equiv g(y^+_{i}, z^+_{i}, p, t_{i})$.
 
Introducing adjoint variables $\lambda$ and $\mu$ (of appropriate dimensions) and using Eq.~\ref{e:FSA_mem}, we have:
\begin{equation}
\begin{split}
\frac{dG}{dp} &= \sum_{i=0}^N \left[
\intgr \left( g_y s + g_z w + g_p \right) dt +
g^-_{i+1} \dtdp{i+1} -
g^+_{i} \dtdp{i} \right. \\
&+ \intgr \lambda^T \left( -\dot s + f_y s + f_z w + f_{y^*} s^+_{i} + f_{z^*} w^+_{i} + f_p \right) dt \\
&+ \left.\intgr \mu^T \left( k_y s + k_z w + k_{y^*} s^+_{i} + k_{z^*} w^+_{i} + k_p \right) dt \right]
\,,
\end{split}
\end{equation}
where we have taken into account that, in mode $\Si$, $s^* = s^+_{i}$ and $w^* = w^+_{i}$.
Next, use integration by parts and rearrange terms to obtain
\begin{equation}
\begin{split}
\frac{dG}{dp} &= \sum_{i=0}^N 
\left[
\intgr \left( g_y + \dot \lambda^T + \lambda^T f_y + \mu^T k_y \right) s\,dt +
\intgr \left( g_z + \lambda^T f_z + \mu^T k_z \right) w\,dt 
\right. \\
&+ \left.
g^-_{i+1} \dtdp{i+1} - g^+_{i} \dtdp{i}
- \left.\left(\lambda^T s\right)\right|_{t_{i}}^{t_{i+1}}
+ \Ii s^+_{i} + \Ji w^+_{i} + \Ki
\right] \,,
\end{split}
\end{equation}
where we define
\begin{align}
\Ii &= \intgr \left(\lambda^T f_{y^*} + \mu^T k_{y^*}\right) \\
\Ji &= \intgr \left(\lambda^T f_{z^*} + \mu^T k_{z^*}\right) \\
\Ki &= \intgr \left(g_p + \lambda^T f_p + \mu^T k_p\right) \, .
\end{align}
As before, the goal is to impose conditions on the adjoint variables $\lambda$ and $\mu$ so that the gradient of $G$ can be computed without the need to compute forward state sensitivities. We first impose that
\begin{subequations}\label{e:Adjoint_mem}
\begin{align}
\dot \lambda^T + \lambda^T f_y + \mu^T k_y + g_y &= 0 \label{e:Adjoint_mem_diff} \\
\lambda^T f_z + \mu^T k_z + g_z &= 0 \label{e:Adjoint_mem_alg} \,,
\end{align}
\end{subequations}
in each interval $[t_i, t_{i+1}]$.
Then, changing summation indices as appropriate and using the fact that the integration limits $t_0$ and $t_f$ are constant (i.e., $\dtdp{0} = \dtdp{N+1} = 0$),
\begin{equation}\label{e:dGdp3}
\begin{split}
\frac{dG}{dp} &= \sum_{i=0}^{N} \Ki + \\
&-\left(\lambda^-_{N+1}\right)^T s^-_{N+1} + \left(\lambda^+_{0}\right)^T s^+_{0} +
I_0^1 s^+_0 + J_0^1 w^+_0 \\
&+ \sum_{i=1}^N \left[
\left( g^-_i - g^+_i\right) \tau_i
-\left(\lambda^-_{i}\right)^T s^-_{i} + \left(\lambda^+_{i}\right)^T s^+_{i} 
+ \Ii s^+_i + \Ji w^+_i
\right]
\end{split}
\end{equation} 
The last term in the above equation dictates the choice of transfer conditions for the adjoint variable $\lambda$ such that all terms involving state sensitivities cancel out.
To do this, we use the conditions of Eq.~\ref{e:FSA_trans_fct_mem} in conjunction with the algebraic equations of Eq.~\ref{e:FSA_mem_alg} evaluated to the left and to the right of the transition at $t_{i}$:
\begin{subequations}
\begin{align}
0 &= \ki{i-1}_y(t^-_{i}) s^-_{i} + \ki{i-1}_z(t^-_{i}) w^-_{i} 
+ \ki{i-1}_{y^*}(t^-_{i}) s^+_{i-1} + \ki{i-1}_{z^*}(t^-_{i}) w^+_{i-1}
+ \ki{i-1}_p(t^-_{i}) \\
0 &= \ki{i}_y(t^+_{i}) s^+_{i} + \ki{i}_z(t^+_{i}) w^+_{i} 
+ \ki{i}_{y^*}(t^+_{i}) s^+_{i} + \ki{i}_{z^*}(t^+_{i}) w^+_{i}
+ \ki{i}_p(t^+_{i}) \,.
\end{align}
\end{subequations}
Also taking into account the dependency of $\dtdp{i}$ on the sensitivities $s^-_i$ (Eq.~\ref{e:FSA_trans_time_mem}), we can formally find a recursion formula for $w^+_{i}$ in terms of the sensitivities at all transition times to the left of $t_i$.

Setting $\lambda(t_f) \equiv \lambda(t^-_{N+1}) = 0$ and computing $\mu(t_f)$ to be consistent with Eq.~\ref{e:Adjoint_mem_alg}, we propagate backward in time the adjoint DAE of Eq.~\ref{e:Adjoint_mem} until the transition time $t^+_N$. At this point, the adjoint variables $\lambda^-_N$, past the transition, are selected as to cancel all terms involving $s^+_N$, $s^-_N$, $w^+_N$, and $w^-_N$ in the last term of Eq.~\ref{e:dGdp3}. This process is repeated in turn for all modes $\Si$ with $i=N,\ldots,0$.

To illustrate this process, consider the following simplifying assumptions:
\begin{itemize}
\item State and state sensitivities are $\mathcal{C}^0$ continuous at all transitions:
    \begin{equation*}
    \begin{split}
    &y^+_i = y^-_i \,, z^+_i = z^-_i \\
    &s^+_i = s^-_i \,, w^+_i = w^-_i
    \end{split}
    \end{equation*}
    Note that this also implies that
    $$
    g^+_i = g^-_i
    $$
\item Initial conditions do not depend on the problem parameters:
    $$
    s_0 = s(t_0) = 0 \,, w_0 = w(t_0) = 0
    $$
\item Without loss of generality (given that $k_z$ is nonsingular), the algebraic equations~\ref{e:FSA_mem_alg} can be written as:
    $$
    w = k_y s + k_{y^*} s^* + k_{z^*} w^* + k_p
    $$
\end{itemize}
Under these assumptions, the gradient of $G$ becomes
\begin{equation*}
\frac{dG}{dp} = \sum_{i=0}^N \Ki 
+ \sum_{i=1}^N \left[
-\left( \lambda^-_i \right)^T s_i + \left( \lambda^+_i \right)^T s_i
+ \Ii s_i + \Ji w_i
\right]
\end{equation*}
The transition conditions for adjoint variables, moving backward in time, are then
\begin{equation*}
\begin{split}
\left(\lambda^-_N\right)^T =&
\left(\lambda^+_N\right)^T + I_{N}^{N+1} + 
J_N^{N+1} k^{(N)}_y \\
\left(\lambda^-_{N-1}\right)^T =&
\left(\lambda^+_{N-1}\right)^T + I_{N-1}^{N} +
\left(J_{N-1}^N + J_N^{N+1} k^{(N)}_{z^*}\right) k^{(N-1)}_y +
J_N^{N+1} k^{(N)}_{z^*} k^{(N-1)}_{y^*} \\
&\vdots \\
\left(\lambda^-_1\right)^T =& \cdots
\end{split}
\end{equation*}
where all partial derivatives are evaluated at the end of the time interval for the corresponding mode; e.g. $k^{(N)}_y \equiv k^{(N)}_y(t^-_N)$.
Note that, even though both the states and sensitivities are continuous across transitions, the adjoint variables are discontinuous with discrete jumps (see also the example in \S\ref{ss:example_EM_memory}).  The complete procedure, also including the required updates to $dG/dp$ at each transition point, is listed in Algorithm~\ref{alg:ASA_mem}.

\begin{algorithm}
\caption{ASA for a Hessenberg index-1 system with memory}\label{alg:ASA_mem}
\begin{algorithmic}[1]
\State Set final conditions: $\lambda = \lambda(t_f) \leftarrow 0$
\State Compute $\mu$ consistent with algebraic constraints $\ki{N}(t_f)$
\State Initialize gradient: $dG/dp = 0 \in \R^{1 \times N_p}$
\State Initialize accumulators: $A \leftarrow 0$, $B \leftarrow 0$
\State Set mode to $\mathbb{S}_{N}^{N+1}$
\For{$i = N,\ldots,0$}
  \State Integrate adjoint system and quadratures over $[t_{i+1}, t_{i}]$:
  \State $\qquad \lambda^+_i, \mu^+_i$ and $\Ii, \Ji, \Ki$
  \State Update: $dG/dp \leftarrow dG/dp + \Ki$
  \If{i = 0}
    \State \texttt{return}
  \EndIf
  \LineComment{{\em Calculate jumps in adjoint variables}}
  \State Set mode to $\mathbb{S}_{i-1}^{i}$
  \State Evaluate derivatives of the algebraic equations: $k_y, k_{y^*}, k_{z^*}, k_p$
  \State Transfer adjoint variables: $\lambda \gets \lambda + \Ii + (\Ji + A) k_y + B$
  \State Compute $\mu$ consistent with algebraic constraints $\ki{i-1}(t^-_i)$
  \State Update: $dG/dp \gets dG/dp + (\Ji + A)k_p$
  \State Update accumulators:
  \State $\qquad B \gets (\Ji + A) k_{y^*}$
  \State $\qquad A \gets (\Ji + A) k_{z^*}$
\EndFor
\end{algorithmic}
\end{algorithm}


\section{Numerical Examples}\label{s:examples}

\subsection{A Simple Hybrid ODE}\label{ss:example_ODE_hybrid}

The first example is reproduced from~\cite{galan1999parametric, tolsma2002hidden} and is a single ODE which transitions between two modes:
\begin{equation}
    \dot x = \left\{\begin{matrix}
    4-x &:& x^3 - 5x^2 + 7x \le p \,, \\
    10-2x &:& \text{otherwise} \,,
    \end{matrix}\right.
\end{equation}
with initial condition $x(t_0) = 0$ and a nominal value $p=2.9$.  The transition function for this example is simply $x^+ = x^-$.

This example is used here to compare the solution and (forward) sensitivity results with those presented in~\cite{galan1999parametric, tolsma2002hidden} and then further compare results of FSA and ASA in computing the derivative of $G=\int_{t_0}^{t_f} x dt$ with respect to the parameter $p$.

\begin{figure}
    \centering
    \includegraphics[width = 0.4\columnwidth]{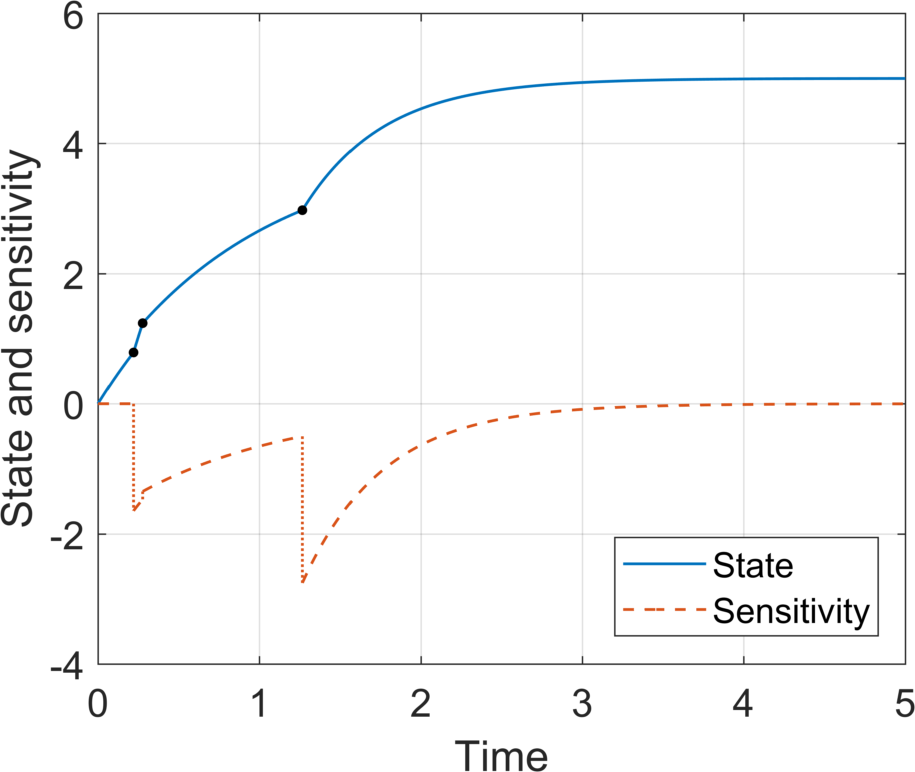}
    \caption{State and sensitivity trajectories for the simple hybrid ODE example.}
    \label{f:example_ODE_FSA}
\end{figure}

Figure~\ref{f:example_ODE_FSA} shows the state variable $x$ and its sensitivity $s = dx/dp$ over the interval $t \in [0, 5]$. These results clearly indicate that explicitly accounting for the discrete jumps in the state sensitivity at the transition times is crucial in obtaining the correct result\footnote{As shown in~\cite{tolsma2002hidden}, ignoring the hidden discontinuities results in $s \equiv 0$ at all times}. The sensitivity discontinuities can be obtained from the differentiated transition function evaluated at a transition time $t^*$:
\begin{equation}
    s^+ + \dot x^+ \frac{dt^*}{dp} = s^- + \dot x^- \frac{dt^*}{dp} \,,
\end{equation}
where the sensitivity of the transition time is obtained from the transition condition $h \equiv x^3 - 5x^2 + 7x - p$:
\begin{equation}
   h_x s + h_x \dot x \frac{dt^*}{dp} + h_p = 0 \,, 
\end{equation}
where all quantities are evaluated right before the transition at $t^*$.
With the state sensitivities such computed, we obtain $\frac{d}{dp}\int_0^5 x dt = -2.31195$.

\begin{figure}
    \centering
    \includegraphics[width = 0.4\columnwidth]{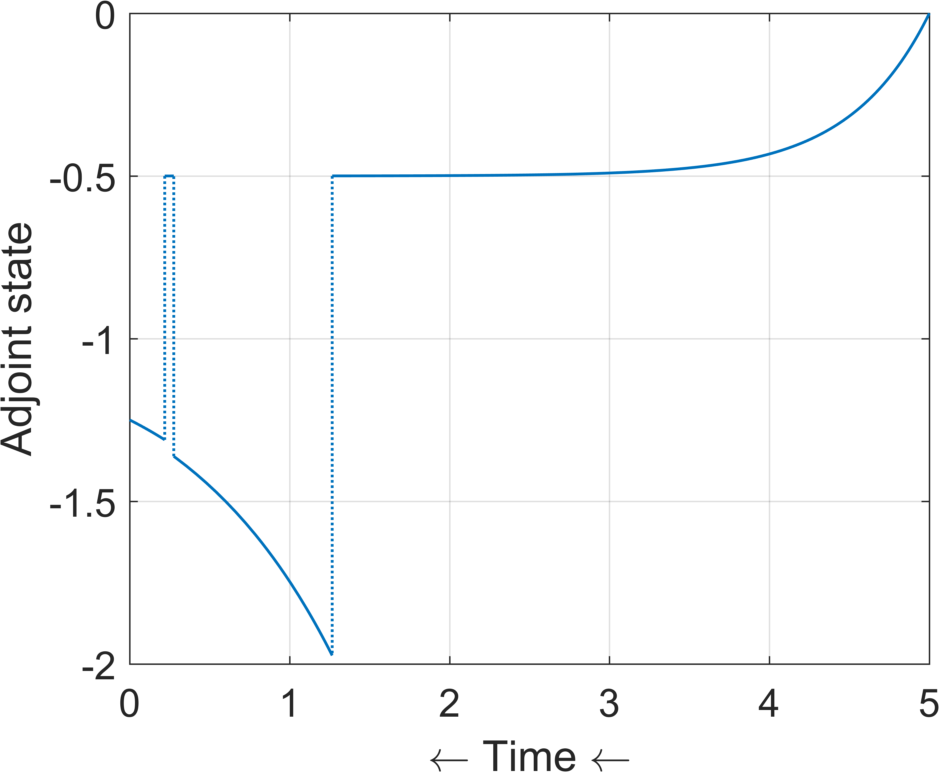}
    \caption{Adjoint variable trajectory for the simple hybrid ODE example.}
    \label{f:example_ODE_ASA}
\end{figure}

Using the adjoint method, the derivative of $G$ is obtained as:
\begin{equation}
    \frac{dG}{dp} = \sum_{i=1}^N \lambda^+_i \left( \dot x^-_i - \dot x^+_i\right) \frac{h_p(t^-_i)}{h_x(t^-_i) \dot x^-_i}
\end{equation}
where the adjoint equations are
\begin{equation}
    \dot \lambda = \left\{\begin{matrix}
    \lambda + 1 &:& x^3 - 5x^2 + 7x \le p \,, \\
    2 \lambda + 1 &:& \text{otherwise} \,,
    \end{matrix}\right.
\end{equation}
with final conditions $\lambda(t_f) = 0$ and transition relationships
\begin{equation}
    \lambda^-_i = \lambda^+_i \left(
    1 - \frac{1}{\dot x^-_i} \left( \dot x^-_i - \dot x^+_i\right)
    \right)\,.
\end{equation}
The trajectory of the adjoint variable is shown in Fig.~\ref{f:example_ODE_ASA}. Note that, for the purpose of evaluating $\frac{dG}{dp}$ backward integration for $\lambda$ could have been stopped once the first transition time $t_1$ is reached. The desired derivative is obtained as $\frac{d}{dp}\int_0^5 x dt = -2.31195$.

\subsection{A Computational Mechanics Problem}\label{ss:example_EM_memory}
This subsection develops a one degree of freedom model using a general constitutive equation with application to computational mechanics. The internal forces are adapted from the {\em exponential model} (EM) of Vaiana et al.~\cite{vaiana2018class}. This material description is intended to fit the rate-independent hysteretic behavior of materials. The EM does not require the use of return-mapping algorithms and can describe a wide range of behaviors.

\begin{figure}
	\centering
	\includegraphics[width = 0.45\columnwidth]{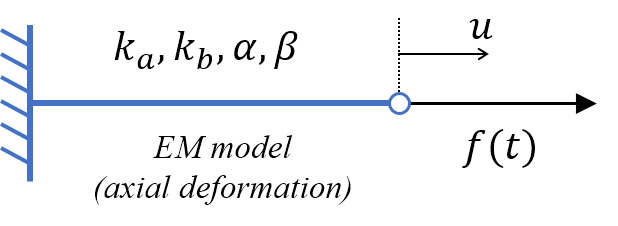}
	\caption{Forced oscillator using the {\em exponential model} (EM) material law.} 
	\label{f:EMBar}
\end{figure}
The dynamics of a single degree of freedom oscillator (see Fig.~\ref{f:EMBar}) can be cast as a Hessenberg index-1 hybrid system with memory:
\begin{subequations} \label{e:EM1}
\begin{align}
\dot u &= v\\ 
m\dot v &= -A z + f(t) \\
z &= \sigma(u, u^*, z^*, p) \,,
\end{align}
\end{subequations}
where $u$ is the deformation displacement, $z$ is the internal stress, $A$ is the cross section area, and $m$ is the mass of the system, concentrated in one node.
The stress is determined by~\cite{vaiana2018class}
\begin{equation}\label{e:EM1_e}
\sigma = -2\beta u + 2 \sinh{\beta u} + k_b u 
- \xi\frac{k_a - k_b}{\alpha} [ e^{-\alpha(u \xi - u^*_i \xi + 2 u_0)} - e^{-2\alpha u_0}] + \xi \bar f
\end{equation}
and depends on the material parameters $p=[k_a, k_b, \alpha, \beta]$.
Transitions are controlled by zero-crossings of the velocity $v$ (since $\xi = sign(v)$) and the dynamics in each mode depend on the state at the last reversal $i$ through the derived memory variable $u^*_i$:
\begin{equation*}
u^*_i = u^* + 2\xi u_0 
+ \frac{\xi}{\alpha}
\ln\left[
  \frac{\xi\alpha}{k_a - k_b}
  \left(
  -2\beta u^* + \sinh{\beta u^*} + k_b u^* + \frac{k_a - k_b}{\alpha} \xi e^{-2\alpha u_0} + \xi \bar f - z^*
  \right)
\right] \,.
\end{equation*}
Here, the quantities $u_0$ and $\bar f$ are (constant) internal model parameters
\begin{equation*}
\begin{split}
u_0 &= -\frac{1}{2\alpha}\ln{\frac{\delta}{k_a - k_b}} \\
\bar f &= \frac{k_a - k_b}{2\alpha} \left(1-e^{-2\alpha u_0}\right)
\end{split}
\end{equation*}
with $\delta = 10^{-20}$.

At each transition point, the state is continuous; i.e., the transition functions are simply $u^+ = u^-$ and $v^+ = v^-$. In each mode, the solution depends on the state at the last transition. Therefore, consistency at each transition point implies that $z^+ = z^-$ and $\dot u^+ = \dot u^-$, $\dot v^+ = \dot v^-$.

\begin{figure}
	\centering
	\includegraphics[width=\textwidth]{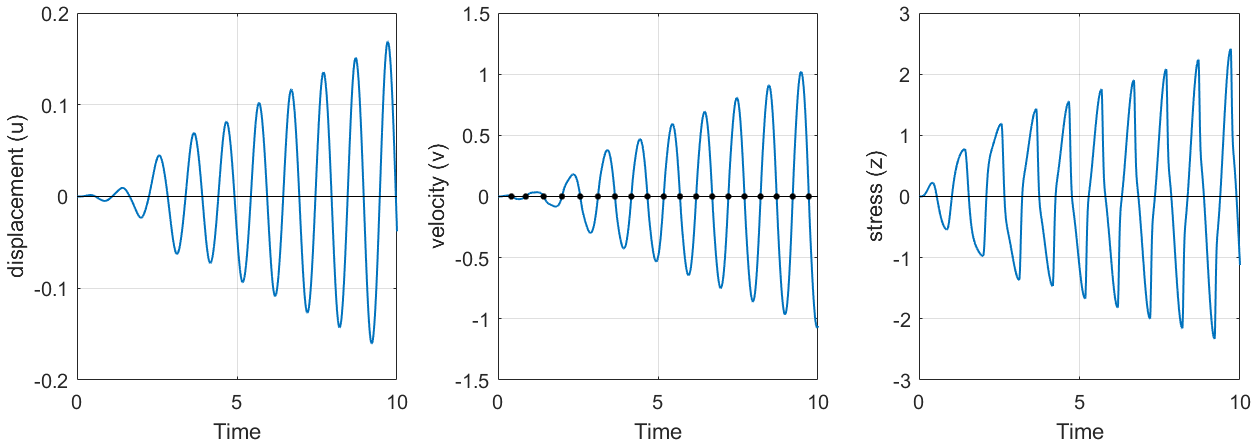}
	\caption{State trajectories for the EM model.}
	\label{f:EM_sim}
\end{figure}

\begin{figure}
    \centering
    \includegraphics[width=0.4\textwidth]{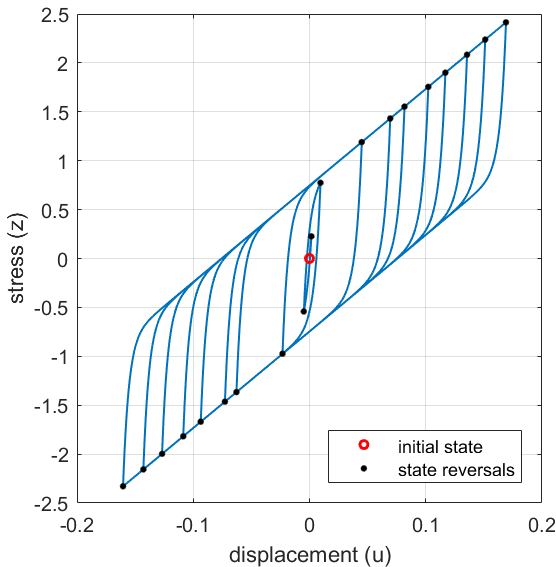}
    \caption{Hysteresis loops and state reversals for the EM model.}
    \label{f:EM_hyst}
\end{figure}

Starting with consistent initial conditions $u(t_0) = v(t_0) = z(t_0) = 0$, with a sinusoidal load  $f(t) = 0.5 t \sin{2\pi t}$ applied to the free end of the bar, the system undergoes multiple reversals over a time interval $t \in [0,10]$ as illustrated in Figs.~\ref{f:EM_sim} and~\ref{f:EM_hyst}. The nominal values for the model parameters were $p=[k_a, k_b, \alpha, \beta] = [32\pi^{2}, \pi^{2}, 205, 0] $, with the mass of the node receiving the excitation $m = 1$ and the cross section area $A=1$. 

\begin{figure}
	\centering
	\includegraphics[width = \columnwidth]{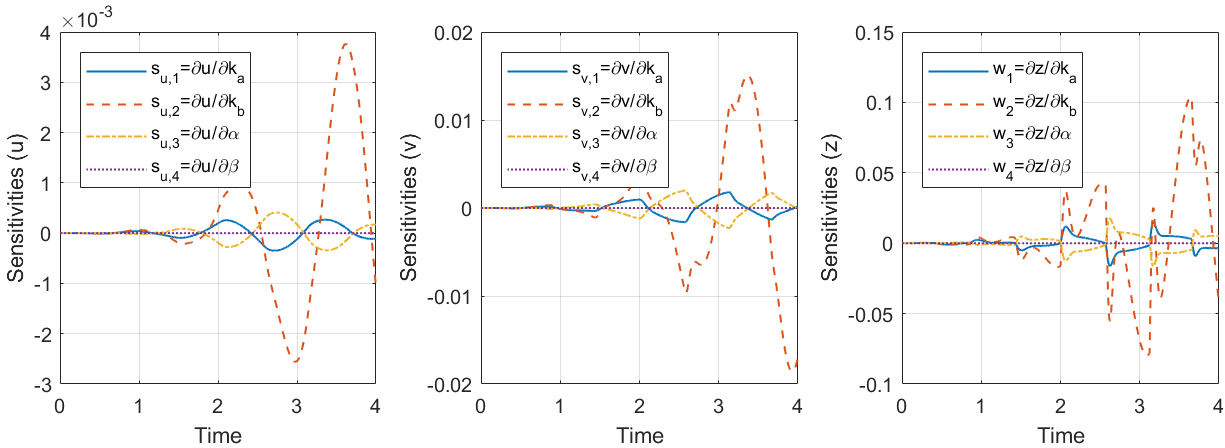}
	\caption{Sensitivity trajectories for the EM model. For clarity, the plots are restricted to the time interval $t\in[0,4]$.}
	\label{f:EM_FSA}
\end{figure}

Following the derivation in \S\ref{s:SA_memory}, the sensitivity equations for this problem take the form
\begin{subequations}
	\begin{align}
	\dot s_u &= s_v \\ 
	m \dot s_v &= -A w \\
	w &= \sigma_u s_u + \sigma_{u^*} s^*_u + \sigma_{z^*} w^* + \sigma_p
	\end{align}
\end{subequations}
with $s_u = \partial u / \partial p \in \R^{1 \times 4}$, $s_v = \partial v / \partial p \in \R^{1 \times 4}$, and $w = \partial z / \partial p \in \R^{1 \times 4}$.  The partial derivatives of the stress function $\sigma$ in the sensitivity algebraic equation were obtained analytically. It is easy to see that, at each transition time $t_i$, the sensitivities of the differential variables are $\mathcal{C}^1$ continuous, while the sensitivity of the algebraic variable is $\mathcal{C}^0$ continuous. See Fig.~\ref{f:EM_FSA}.

As output functional, we consider a measure of total deformation over the time interval considered defined as $G(p)=\int_{t_0}^{t_f} u^2 dt$. During integration of the DAE and sensitivity equations, additional quadratures are carried out to compute $G$ and its gradient with respect to the problem parameters $p$. 
With the current choice of nominal model parameter values and a final simulation time $t_f = 10$, the computed value is $G=0.04994$.

\begin{figure}
	\centering
	\includegraphics[width=0.5\textwidth]{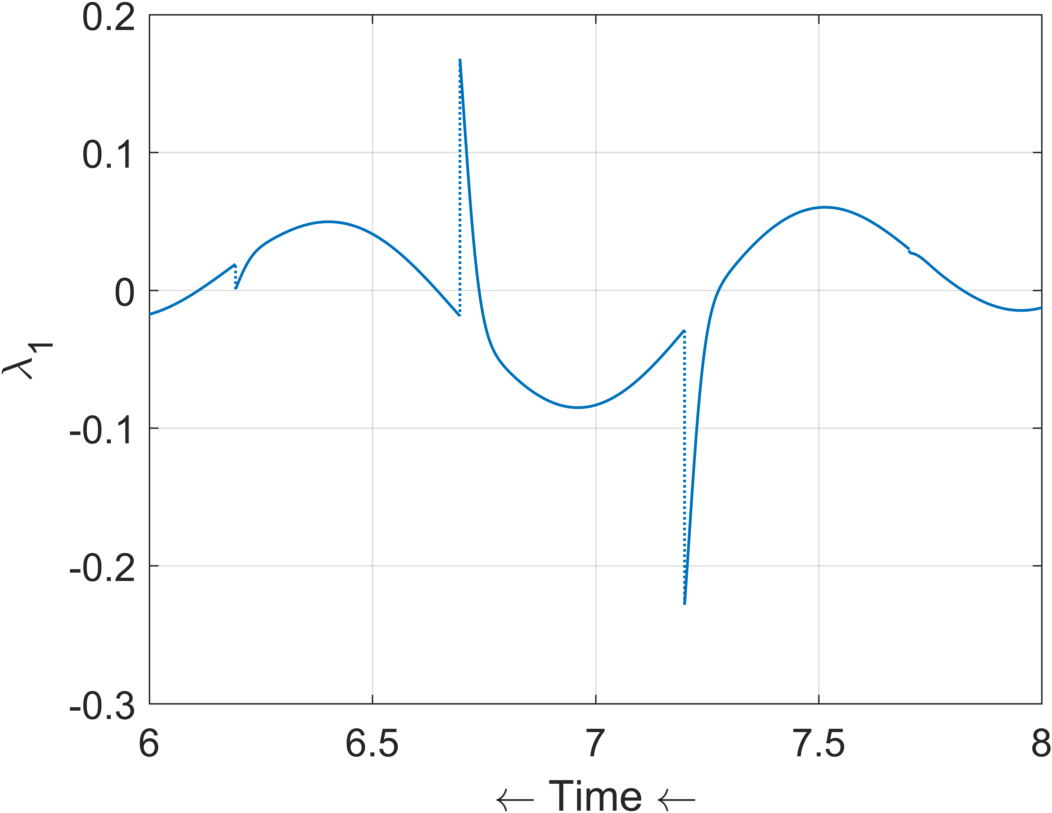}
	\caption{Trajectory of $\lambda_1$ with discontinuities at transition points. For clarity, the plot is restricted to the time interval $t\in[6,8]$.}
	\label{f:EM_adj}
\end{figure}

Finally, the procedure of \S\ref{ss:ASA_memory} described in Algorithm~\ref{alg:ASA_mem} was employed to compute the gradient $dG/dp$ with the ASA approach.  The adjoint systems for this problem assume the form:
\begin{subequations}
	\begin{align}
	\dot \lambda_1 &= \sigma_u \mu - 2 u \\ 
	\dot \lambda_2 &= -\lambda_1 - \sigma_v \mu \\
	0 &= \mu - \lambda_2 / m
	\end{align}
\end{subequations}
For this example problem, the recurrence relationship controlling the transfer of adjoint variables at transition points has the form:
$$
w^i = w^{i-1} + \sigma_u \left( s_u^i - s_u^{i-1}\right) + \sigma_p \,.
$$
Starting with the consistent conditions $\lambda(t_f) = [0 \,,\, 0]$ and $\mu(t_f) = 0$ at $t_f = 10$, the adjoint solution is shown in Fig.~\ref{f:EM_adj}. It is important to highlight that, even though both the states and sensitivity trajectories in this problem are $\mathcal{C}^0$ continuous across transitions, the adjoint variables are discontinuous with discrete jumps which must be properly accounted for and incorporated during the backward-in-time propagation of the associated hybrid adjoint system.

As an additional check on the output functional gradient computed using FSA and ASA approaches, we also provide below an estimate based on (forward) finite differences, using a relative perturbation factor $\bar\epsilon = 10^{-4}$.
The gradient of $dG/dp$ obtained with different techniques is:
\begin{align*}
    \text{FD:}\,&  [-1.338\cdot10^{-5} , 3.267\cdot10^{-3} , -1.534\cdot10^{-6} , -6.07\cdot10^{-9}] \\
    \text{FSA:}\,& [-1.337\cdot10^{-5} , 3.266\cdot10^{-3} , -1.518\cdot10^{-6} , 0] \\
    \text{ASA:}\,& [-1.335\cdot10^{-5} , 3.267\cdot10^{-3} , -1.540\cdot10^{-6} , 0] \,,
\end{align*}
consistent with the integration tolerances and the relative finite difference perturbation factor. We note that integration of all DAEs (both forward and backward in time) was conducted with the {\tt ode15s} Matlab solver, using a relative tolerance of $10^{-8}$ and an absolute tolerance of $10^{-12}$. 


\section{Conclusions}

We have derived the sensitivity equations for forward sensitivity analysis and the adjoint equations for adjoint sensitivity analysis for hybrid discrete -- continuous systems whose dynamics between transition points are described by index-0, index-1, and index-2 Hessenberg DAEs. Using a general framework for describing such systems, we provide sufficient solvability conditions and derive appropriate boundary (final) conditions, as well as the transfer relationships at transition points for the adjoint variables.

Furthermore, we have extended the derivation to include so-called hybrid systems with memory, where the dynamics depend explicitly on the state at the last transition and demonstrated this approach on a computational mechanics problem. Ongoing follow-up work focuses on extensions to more complex dependencies on past history, with application to different material constitutive laws, such as the Menegotto-Pinto model.




\newpage
\bibliographystyle{unsrt}
\bibliography{sa}


\end{document}